\magnification=\magstep1
\centerline {\bf Chern classes and the periods of mirrors.}
\bigskip 
\centerline {{by Anatoly Libgober}%
\footnote{${}^*$}{
Supported in part by NSF}}
\bigskip 
\centerline {Department of Mathematics}
\centerline {University of Illinois at Chicago}
\centerline {851 S.Morgan, Chicago, Ill, 60607}
\centerline {{\it e-mail:} libgober@math.uic.edu}
\bigskip {\bf Abstract.} {\it We show how Chern clases of a Calabi
Yau hypersurface in a toric Fano manifold can be expressed in terms 
of the holomorphic at a maximal degeneracy point period of its mirror.
We also consider the relation between the Chern classes and the periods 
of mirrors for complete intersections in Grassmanian Gr(2,5).} 
\bigskip 
{\bf 0.Introduction.} 
\smallskip
A relation between Chern classes of 3-dimensional Calabi Yau
complete intersections in a product of weighted projective spaces 
and periods of holomorphic 3-forms on the mirror was described in [HKTY] 
where the following  ``remarkable identities'' were 
observed. Let $X$ be a 3-dimensional Calabi-Yau complete intersection  
in a product of $k$ weighted projective spaces. Let 
$ w_0(z_1,..,z_k)=\Sigma_{n_i \ge 0}
c(n_1,...n_k)z_1^{n_1}...z_k^{n_k}$
be the period of a holomorphic 3-form on  $X$ 
which admits a holomorphic extension into the point
with a maximally unipotent monodromy normalized so that it
takes there the value 1. Then 
$c(n_1,..,n_k)$ is a product of $(l_j(n_1,..,n_k)!)^{\pm 1}$ where $l_j$ are
linear forms with integer coefficients. Define $c(...,\rho_i,..)$
replacing in $c(n_1,...,n_k)$ 
each $(l_j)!$ by $\Gamma(l_j(..,\rho_i,...)+1)$.   
For $i=1,..,k$ let $J_i$ be the class of the 
 pull back on $X$ of the Kahler form of $i$-th projective space which
product contains $X$. Let $K_{ijk}$ be the Yukawa coupling,
$\partial_{\rho_i}={1 \over {2 \pi i}} {\partial \over {\partial
\rho_i}}, D^{(2)}_i={1 \over 2}
K_{ijk}\partial_{\rho_j}\partial_{\rho_k}, 
D^{(3)}=-{1 \over 6}K_{ijk}\partial_{\rho_i}  \partial_{\rho_j}
 \partial_{\rho_k}$.
Then (cf. (4.21) in [HKTY])
$$\int_X c_2 \wedge J_i=-24 D^{(2)}c(\rho_1,...,\rho_k) 
\vert_{(0,..,0)} \eqno (0.1)$$
$$\int_Xc_3= i {{2 \pi}^3 \over {\zeta (3)}} D^{(3)}c(\rho_1,...,\rho_k)
\vert_{(0,..,0)} \eqno (0.2)$$ 
 The purpose of this note is to  extend this relationship 
to arbitrary dimension. We shall consider the case of 
hypersurfaces in toric Fano manifolds.
For these Calabi Yau manifolds, 
 we show that certain linear combinations of 
Chern classes can be expressed in terms of the period of the mirror
holomorphic at the maximal degeneracy point. These combinations
form a Hirzebruch's multiplicative sequence i.e. obtained by applying 
Hirzebruch construction [Hi] to the series ${1 \over {\Gamma(1+z)}}$. 
This multiplicative $\Gamma$-sequence,
in turn, determines the Chern classes of the manifold (cf. lemma below). 
\par More precisely we prove the following.
Let $\Delta$ be a reflexive Fano polytope
(cf. [B1],[B2]) of dimension $d$ (in other words 
 $\Delta$ is a simplicial polytope 
with: (a) vertices belonging to an integral lattice $M$ 
of rank $d$; (b) the origin being the only intersection point
of the interior of $\Delta$ and the lattice $M$; (c) such that for 
each $(d-1)$-dimensional face its vertices
form a basis of $M$ and (d) the equation of this face is $l=-1$ 
where $l$ is a linear form on $M$;). 
The corresponding toric variety $X_{\Delta}$ is a Fano manifold i.e. 
its first Chern class is ample and the hypersurface which is the  
zero locus of a section of the line bundle on $X_{\Delta}$  
corresponding to $\Delta$ is a Calabi Yau hypersurfaces $V_{\Delta}$.
Let $\Delta^* \subset M^* \otimes {\bf R}$ ($M^*=Hom(M,{\bf Z})$) 
be the polar polytope and
 $f_{\Delta^*}$ be a generic linear combination of characters of 
$M \otimes {\bf C}^*$ which belong to $\Delta^*$. 
Then a period of the affine hypersurface $f_{\Delta^*}=0$ in 
$Hom (M^*, {\bf C}^*)$  i.e.  
$\int_{\gamma} {1 \over {f_{\Delta^*}}} \prod {{dx_i} \over
{x_i}}$ (where $x_i$'s are coordinates in $Hom (M^*,{\bf C}^*)$
corresponding to a choice of a basis in $M$ 
and $\gamma$ is a $d$-cycle in the complement to $f_{\Delta^*}=0$
in $Hom(M^*,{\bf C}^*)$)    
satisfies a system of Picard Fuchs PDE. In appropriate 
partial compactification of $Hom (M^*,{\bf C}^*)$ 
(constructed in [HLY]) there is a maximal degeneracy point i.e. 
such that this system has only one, up to a constant factor, solution which 
admits a holomorphic extension  into this point. Near the maximal 
degeneracy point
this period has form $\Sigma_{n_1,...n_r \in {\bf N}^{\ge 0}}
 \prod_i l_i(n_1,..,n_r)^{\pm 1} z_1^{n_1} \cdot \cdot
\cdot z^{n_r}$ where $l_i$ are linear forms with integer coefficients. 
Let $c_{\rho_1,...\rho_r}=
{\Gamma}(l_i(\rho_1,..,\rho_r)+1)^{\pm 1}$.
 We shall identify the tangent space to the compactification
of the moduli space of hypersurface $X_{\Delta^*}$ 
at a maximal degeneracy point with $H^2(X_{\Delta}, {\bf C})$ and 
coordinates $(t_1,...,t_r)$ of $H^2(X_{\Delta},{\bf C})$ in a basis 
$J_1,...,J_r \in H^2(X_{\Delta})$ with the coordinates
$(x_1,...,x_r)$ in the partial compactification of the moduli space of 
hypersurfaces $V_{\Delta}^*$ near a maximal degeneracy point.  
Let $K_{i_1,...,i_d}=\int_{X_{\Delta^*}} \Omega \wedge {{\partial^d
\Omega} \over {\partial x_{i_1}...\partial x_{i_d}}}$ be the $d$-point 
function corresponding to 
$V_{\Delta^*}$ (cf. [GMP],[BvS]). Its value at the maximal degeneracy point is normalized
to be equal to $\int_{V_{\Delta}} J_{i_1} \wedge ... \wedge J_{i_d}$
(with such normalization it yields the data of enumerative geometry of
rational curves on $X_{\Delta}$ cf. [GPM]). 
Then we have:
\bigskip {\bf Theorem.} If the assumption (**) below is satisfied, the degree $k$ polynomial 
$ \hfil \break Q_k(c_1,...,c_k)$ in 
the Hirzebruch's multiplicative sequence corresponding to the 
series $1 \over {\Gamma(1+z)}$ of $V_{\Delta}$ 
satisfies: $$\int_{X_{\Delta}} Q_k(c_1,...,c_k) \wedge J_{i_1} \wedge ...
\wedge
J_{i_{n-k}}=$$ $$\sum_{(j_1,..,j_k,i_1,..,i_{n-k})} { 1 \over {k!}}
{{\partial^kc(\rho_1,..,\rho_k)} \over {\partial \rho_{j_1} ... 
\partial \rho_{j_k}}} \vert _{\rho_{j_1}=0,...,\rho_{j_k}=0}
K_{j_1,..,j_k,i_1,...,i_{n-k}} \eqno (0.3)$$
 where the summation is over all permutations $(j_1,....,i_{n})$ of
$(1,..,n)$.
\bigskip As a corollary we obtain an expression for the Chern classes 
in terms of the period and the d-point function of the mirror. 
Multiplicative $\Gamma$-sequence is discussed in section 1, the proof
of the theorem in section 2 where also special cases of it are made
explicit. In section 3 we discuss the simplest example  of the theorem
and relation between the periods on Chern classes for manifolds which 
mirror was constructed recently in [BCKS] via conifold transitions.
\eject
\bigskip {\bf 1. Multiplicative $\Gamma$-sequence}.
\smallskip Let $Q(z)={1 \over {\Gamma (1+z)}}$. 
Since $Q(0)=1$ there is a well deifined  multiplicative sequence
corresponding to this power series (cf. [Hi]). 
The latter associates with a complex manifold $M$ having 
Chern classes $c_i, i=1,...dim \ M$  weighted 
homogeneous polynomials $Q_i (c_1,...c_i)$ 
of degree $2i$ where $deg c_i =2i$. 
\bigskip {\bf Lemma.} The coefficient of $c_i$ in $Q_i$ is equal to 
$\zeta(i)$ for $i \ne 1$ and to the Euler 
constant $\gamma=lim_{n \rightarrow \infty} \Sigma {1\over n}-ln \ n$
for $i=1$. 
In particular it is non zero and hence $c_i$ is a  
the polynomial in $Q_j(c_1,...c_i), (j \le i)$. 
\bigskip {\bf Proof.} The coefficient of 
$c_i$  in $Q_i(c_1,..,c_i)$ is $s_i$ where 
$$1-z{d \over {dz}}log Q(z)=\Sigma_{i=0}^{\infty} (-1)^is_iz^i \eqno (1.1)$$
(cf. [H] sect. 1.4). In the case of $\Gamma$-sequence we have (cf. [E] p. 45) 
$$log \Gamma (1+z)=-\gamma z +\Sigma_{i=2}^{\infty} (-1)^i \zeta(i)
z^i/i \eqno (1.2)$$ This yields the lemma.
\bigskip It isn't hard to work out explicitly the polynomials 
$Q_i$ (using for example the recurrence formulas from [LW]). For small
$i$ we have the following ($\gamma$ 
is the Euler constant $lim_{m \rightarrow \infty} \Sigma_1^m {1\over n}-
ln \ m $): 
   $$Q_1(c_1)=\gamma c_1 \eqno (1.3)$$
    $$Q_2(c_1,c_2)=-{1 \over 2} \zeta(2)c_1^2+\zeta(2)c_2+\gamma^2c_1^2 
 \eqno(1.4)$$
    $$Q_3(c_1,c_2,c_3)=\zeta(3)c_3-(\zeta(3)-\zeta(2)\gamma)c_1c_2+
    ({1 \over 3} \zeta(3)+{1 \over 6} \gamma^3)c_1^3 \eqno (1.5)$$
   $$Q_4(c_1,c_2,c_3,c_4)=\zeta(4)c_4 +{1\over 2}(\zeta(2)^2-\zeta(4))c_2^2
    +(-\zeta(4)+\zeta(3)\gamma)c_1c_3+$$
   $$(\zeta(4)-\zeta(3)\gamma+{1
    \over2}\zeta(2)\gamma^2-{1 \over 2}\zeta(2)^2)c_1^2c_2+
   (-{1 \over 4}\zeta(4)-{1\over 4}\zeta(2)\gamma^2+{1
    \over3}\zeta(3)\gamma+
      {1 \over 8}\zeta(2)^2+ {1 \over 24} \gamma^4)c_1^4 \eqno (1.6)$$
In particular for a Calabi Yau manifold we have the following 
polynomials: $$Q_2(c_2)=\zeta(2)c_2, Q_3(c_2,c_3)=\zeta(3)c_3, 
Q_4(c_2,c_3,c_4)={1 \over 2}(\zeta(2)^2-\zeta(4))c_2^2+\zeta(4)c_4 \eqno
    (1.7)$$   
\bigskip {\bf 2. Periods for the mirrors of hypersurfaces in 
toric varieties}. 
\bigskip \par We shall use the construction of the maximal degeneracy
point and notations from [HLY]. 
Let $M$ be a lattice ${\bf Z}^d$ and $N$ denotes its dual.
A cone in $N\otimes {\bf R}$ is called large if it is generated by 
finitely many vectors from $N$ and its dimension is equal 
to $rk N$. A large cone is called regular if generators of 
its one dimensional faces (edges) form a basis of $N$. 
A fan $\Sigma$ in $N \otimes {\bf R}$ is called regular if 
all its large cones are regular. We consider  
the fan $\Sigma$ corresponding to a Fano polytope 
$\Delta \subset M \otimes {\bf R}$. 
The intersections of the cones from $\Sigma$ and the polytope $\Delta^*$
polar to $\Delta$ form a (maximal, cf. [HLY]) triangulation of the latter. 
As was mentioned the corresponding toric variety $X_{\Delta}$
is a non singular Fano manifold.   
\par Let $\Sigma(1)$ be the set of primitive elements on the 
edges of the fan $\Sigma$. Since $\Delta$  is a Fano polytope
the following property (*) (assumed in [HLY], cf. 4.7) is satisfied:
 $\Sigma (1) \subset \partial conv(\Sigma(1) $. 
We shall denote the points of $\Sigma(1)$ as 
$\mu_1,...,\mu_p$ and put $\mu_0=(0,..,0)$.
\par Let ${\cal N}$ be a lattice endowed with the basis $n_0,...n_p$ 
elements of which are in one to one correspondence with the set 
${\cal A}=\{ \mu_0,...,\mu_p \} \subset {\cal N}$. 
Let $\bar {\cal A}=\{\bar \mu_0,..., \bar \mu_p \}$ be
the collection of points in $\bar { N}={ N} \otimes {\bf R} \oplus {\bf R}$ 
all belonging to the affine
hyperplane $N \otimes {\bf R} \times 1$ and such that the image 
of $\bar \mu_i$ in $N \otimes R$ is $\mu_i$. 
 Denote by  $R: {\cal N} \rightarrow {\bar N}$
the map given by $R(n_i)=\bar \mu_i$ and 
let $L_{\cal A}=Ker R$ be the lattice of relations among $\bar \mu_i$'s
(i.e. $L_{\cal A}=\{(l_0,..l_p) \in {\cal N}
 \vert \Sigma_i l_i\bar \mu_i=0\}$).
Let us put $\bar M=\bar N^*$ and view it as a subspace of the space of functions on 
$\bar {\cal A}$. The quotient ${\cal N}^*/\bar M$ has
the canonical identification
with $L_{\cal A}^*$ (i.e. the dual of $L_{\cal A}$).
 The space ${\cal N}^*$ supports a natural 
fan edges of which corresponds to triangulations of $conv ({\cal A})$ 
with all simplices having 
vertices belonging to ${\cal A}$. The image of this 
fan in $L_{\cal A} \otimes {\bf R}$ is called 
the secondary fan $S\Sigma$ (cf. [GKZ]). Moreover the fan $S\Sigma$ admits 
a natural refinement called the Grobner fan (cf. [HLY] section 4).
We shall assume that:
\bigskip (**) the cone $C({\cal A})$ 
in $L_{\cal A}^*$ generated by 
$\mu_1,..,\mu_p$ is a cone of both the secondary and Grobner fans 
and that the generators of the edges of this cone form a basis of the 
integral lattice of $L_{\cal A}$.  
\bigskip The space $PL(\Sigma)$ of piecewise linear functions on 
$N \otimes {\bf R}$ linear
on each cone of $\Sigma$ can be identified with a subspace of ${\cal
N}^*$.  Moreover $L^*_{\cal A}=
{\cal N}^*/M \otimes {\bf R}=PL(\Sigma)/{\bar M} \otimes {\bf R}$
and the latter is canonically isomorphic to $H^2(X_{\Delta},{\bf R})$.
Lefschetz theorem identifies this group with $H^2(V_{\Delta},{\bf R})$
 if $d \ge 3$. 
In this identification the Kahler cones of $X_{\Delta}$ and $V_{\Delta}$ are
identified with $C({\cal A})$.
\par According to [HLY] a consequence of (**) is that in the 
partial compactification of 
 $\hfil \break {{\bf C}^*}^{\cal A}/{\bf C}^* \times Hom (N,{\bf
C}^*)$, given by the fan consisting of the cones in the closure of
the Kahler 
cone, the point corresponding to the cone $C({\cal A})$ 
is a maximal degeneracy point for the GKZ system:  
 $$(\prod_{l_{\mu}>0} ({\partial \over {\partial a_{\mu}}})^{l_{\mu}}-
 \prod_{l_{\mu}<0} ({\partial \over {\partial a_{\mu}}}))^{-l_{\mu}}
\Pi (a)=0 (l \in  L_{\cal A}),
(\Sigma_{\mu}<u,\bar \mu>a_{\mu}{\partial \over {\partial a_{\mu}}}-
<u,\beta>)\Pi(a)=0 \eqno (2.1)$$ 
where $\beta=(-1,0,...,0)$ and $u \in \bar M\otimes {\bf R}$.
In other words there is only one solution
of (2.1) which admits a holomorphic extension in a neighborhood of the 
point of compactification corresponding to the cone $C({\cal A})$.
Moreover the period: $$\Pi_{\gamma}(a)=\int_{\gamma}{1 \over {f_{\cal
 A}}} \prod {{dX_i} \over {X_i}}, \ \ f_{\cal A}=\Sigma a_{\mu}X^{\mu} \eqno (2.2)$$
($\gamma$, as above,
 is a $d=rk N$ cycle in the complement in $Hom (N,{\bf C}^*)$
to the hypersurface $f_{\cal A}(X,a)=0$) is a solution of the system
 (2.1). Since we assume that $\Delta$ is Fano in fact any solution of (2.1) is 
a period (cf. [H]).  
\par On the other hand  we have the following series representation:
$$a_0 \Pi_{\gamma}(a)=(2\pi i)^d \Sigma_{l_1\mu_1+...l_p\mu_p \ge 0,
l_1,..,l_p = 0} (-1)^{l_1+...l_p}
{{(l_1+...l_p)!} \over {(l_1)!...(l_p)!}} {{a_1^{l_1} \cdot \cdot \cdot
 a_p^{l_p}} \over {a_0^{l_1+....l_p}}} \eqno (2.3) $$
The summation in the latter can be changed into summation 
over the Mori cone in $L_{\cal A}$ by assigning to $(l_1,...,l_p)$
the vector $(l_0,l_1,...,l_p) \in L_{\cal A}$ where $l_0=-l_1-...-l_p$.
If $l^{(1)},...,l^{(p-d)}$ 
is a basis in the Mori cone then in corresponding 
canonical coordinates (cf. (3.1) [HLY])
$x_k=(-1)^{l_0^{(k)}}a^{l^{(k)}}$ we have (cf. (5.12) in [HLY]) 
$$a_0\Pi_{\alpha}(x_1,...,x_{p-d})=$$ $$\Sigma_{m_1,..,m_{p-d}, \Sigma
m_kl^{(k)}_0 \le 0}
{{\Gamma(-\Sigma m_kl^{(k)}_0+1)}\over
{\Gamma(\Sigma m_kl^{(k)}_1+1) \cdot \cdot \cdot 
\Gamma (\Sigma m_kl^{(k)}_p+1)}}x_1^{m_1},...,x_{p-n}^{m_{p-n}}
 \eqno (2.4)$$
Let $J_k (i=1,...p-d)$ be elements of $L^*_{\cal A}$ forming the 
basis dual to $l^{(k)}$ and $D_i, (i=1,..,p)$ be the cohomology classes
in $H^2(X_{\Delta},{\bf Z})$ dual to codimension one orbits of 
${X}_{\Delta}$. Since $D_i$ correspond to the generators of 
one dimensional cones of $\Sigma$ and under identification 
$H^2({V}_{\Delta},{\bf Z})$ correspond to $\mu_i$ 
we have $D_i=\Sigma_k J_kl^{(k)}_i, (i=1,..,p)$ and the total Chern
class of the Calabi Yau hypersurfacs in $X_{\Delta}$, 
which has as the the dual cohomology class  $D_1+...D_p$, is 
the sum of the terms of degree  not exceeding $d$ in the expansion of 
$${{(1+D_1) \cdot \cdot \cdot (1+D_p)} 
\over {(1+D_1+...+D_p)}} \eqno (2.5)$$ 
Since in (2.4) one can view $m_k$ as elements of $L_{\cal A}^*$ 
and hence identified them with $J_k$`s, it follows that  
 the term  $Q_k(c_1,..,c_k)$ of the degree $k$ of
$\Gamma$-sequence coinsides with the term of degree $k$ in the formal series  
$$C(J_1,..J_{n-p})={{\Gamma(-\Sigma J_kl^{(k)}_0+1)} \over 
{\Gamma(\Sigma J_kl^{(k)}_1+1) \cdot \cdot \cdot 
\Gamma (\Sigma J_kl^{(k)}_p+1)}} \eqno (2.6) $$
The latter is equal to $$\sum_{j_1,..,j_k} {1 \over {k!}}
{{\partial^k C(J_1,...,J_k)}
\over {\partial J_{j_1},..., \partial J_{j_k}}}J_{j_1}...J_{j_k} \eqno (2.7)$$
The theorem follows from comparison (2.4) and (2.6). 
\bigskip {\bf Corollary} 1.(cf. (0.1),(0.2),[HKTY].) 
Let $X$ be a Calabi Yau hypersurface of
dimension 3. Then $$\int_X c_2 \wedge J_i = {3 \over {\pi^2}}K_{ijk}
{{\partial^2 c(0,...0)} \over {\partial x_i \partial x_j \partial x_k}}
\eqno (2.8) $$
 $$\int_X c_3 = {6 \over {\zeta(3)}}K_{ijk}
{{\partial^3 c(0,...0)} \over {\partial x_i \partial x_j \partial x_k}}
\eqno (2.9)$$
\par 2. Let $X$ be a Calabi-Yau hypersurface of dimension
4 in a non-singular toric Fano manifold. Then 
  $$\int_X c_2 \wedge J_i \wedge J_j={3 \over {\pi}^2} K_{ijkl} {{\partial c}^2(0,...,0) 
\over {\partial \rho_k \partial \rho_l}},
\int_X c_3 \wedge J_l= {6 \over {\zeta(3)}}K_{ijkl}
{{\partial^3 c(0,...0)} \over {\partial x_i \partial x_j \partial x_k}}
 \eqno (2.10)$$
 $$\int_X ({1 \over 2} {\zeta(2)^2-\zeta(4)}) c_2^2+\zeta(4) c_4=
{1 \over {24}}
K_{ijkl} {{\partial c}^4(0,..,0) \over {\partial \rho_i \partial \rho_j 
\partial \rho_k \partial \rho_k}} \eqno(2.11)$$
\bigskip These identities are a consequences of (0.3) and identities (1.7).  
\bigskip {\bf 3. Concluding remarks.} 
\smallskip 1. {\bf Example.} For a hypersurface $V_d$  of degree $d$ in ${\bf P}^d$
the Chern polynomial is the sum of the terms degree less than  $d+1$ 
in the series ${{(1+h)^{d+1}} \over {(1+dh)}}$ i.e. the $\Gamma$-sequence
for $V_d$ is the sum of the terms of less than $d+1$ in 
${{\Gamma (1+dh)} \over {\Gamma (1+h)^{d+1}}}$. (0.3) is a consequence
of the fact that after replacing $h$ by $m \in {\bf Z}^{+}$ 
the latter becomes  the coefficient of the holomorphic 
at the maximal degeneracy point period of the d-form on the mirror of $V_d$. 
\bigskip 2.Let $X_{1,1,3}$ (resp. $X_{1.2.2}$ 
be a non singular complete intersection of the Grassmanian
$Gr(2,5)$ embedded via Plucker embedding in ${\bf P}^9$ and  
 hypersurfaces of degrees $(1,1,3)$ (resp. $(1,2,2)$). 
It is shown in  [BKCS] that holomorphic period of the mirror has  
presentation
 $$\sum_m \left [{(3m)!} {m!}^2 \sum_{r,s} {1 \over {m!}^5} {m \choose r} 
{s \choose r} {m \choose s}^2 \right ] z^m  \eqno (3.1)$$
resp.
 $$\sum_m \left [{m!}{(2m)!}^2 \sum_{r,s} {1 \over {m!}^5} {m \choose r} {s \choose 
r} {m \choose s}^2 \right ] z^m \eqno (3.2)$$  
Though the corresponding differential equations are not of hypergeometric
type (since they have more than three singular points) and the ratio of 
coefficients is not a rational function of $m$, nevertheless
the ratio of the coefficients of the series (3.1) and (3.2) 
is equal to the value of the ratio of Hizrebruch $\Gamma$ sequences for 
$X_{1,1,3}$ and $X_{1,2,2}$. Indeed the Chern polynomials
of the these manifolds are restrictions on 
cohomology of corresponding manifolds of 
${{c(Gr(2,5)} \over {(1+h)^2(1+3h)}}$
and ${{c(Gr(2,5)} \over {(1+h)(1+2h)^2}}$ respectively where 
$c(Gr(2,5))$ is the Chern polynomial of  the Grassmanian and $h$ is the 
cohomology class of the hyperplane section i.e. the ratio of the
$\Gamma$-sequences is ${{\Gamma(3h+1) \Gamma(h+1)} \over {\Gamma(2h+1)}}$.  
which for after replacing $h$ by $m \in {\bf Z}^{+}$  is equal to the 
ration of the coefficients of the periods of mirrors of corresponding
Calabi Yau manifolds
\eject
\bigskip  
\bigskip 
\centerline {\bf References}
\bigskip [B1] V.Batyrev, Dual polyhedra and mirror symmetry for 
Calabi Yau hypersurfaces in toric vareties, 
J. of Alg. Geom., 3 (1994), 493-535.
\bigskip [B2] V.Batyrev, On the classification of toric Fano 4-folds,
math.AG 9801107.
\bigskip [BKCS] V.Batyrev, I.Ciocan-Fontanine, B. Kim, D. v. Straten,
Conifold Transitions and mirror symmetry for Calabi Yau complete
 intersections in Grassmanians, alg-geom/9710022.
\bigskip [BvS] X.Batyrev, D. v. Straten, Generalized hypergeometric functions and rational 
curves on Calabi Yau complete intersections in toric varieties.
alg-geom 9307010.
\bigskip [E] A.Erdelyi ed., Higher transcendental functions, vol. I.
McGraw Hill Book Co. New York, 1953. 
\bigskip [GMP] B.Greene, D.Morrison, R.Plesser, Mirror manifolds in
higher dimensions. 
\par Comm. Math. Phys.,173 (1995),1699-1702.
\bigskip [Hi] F.Hirzebruch, Topological methods in algebraic geometry, 
Springe-Verlag, 1966.
\bigskip [H] S.Hosono, GKZ systems, Grobner Fans and moduli spaces 
of Calabi Yau hypersurfaces, alg-geom/ 9707003.  
\bigskip [HKTY] S.Hosono, A.Klemm, S.Theisen, S.-T. Yau. Mirror
symmetry, mirror map and applications to complete intersections
Calabi Yau spaces.
\bigskip [HLY] S.Hosono, B.H.Lian and S.-T. Yau. Maximal degeneracy
points of GKZ systems, J. Amer. Math. Soc, vol.10, No.2 (1997), 
p. 427-443. 
\bigskip [GKZ] I.M.Gelfand, M.Kapranov, A.Zelevinski, Discriminants,
Resultants and Multidimensional Determinants, 
Birhauser, Boston, 1994.
\bigskip [LW] A.Libgober, J.Wood, Uniqueness of the complex structure on
Kahler manifolds of certain homotopy type, J. Diff. Geom, 32 (1990), p. 139-154.
\end